\documentclass[11pt]{article}
\usepackage{amsmath}
\usepackage{amssymb}
\usepackage{amscd}
\usepackage{amsthm}

\def\endproof{\relax\ifmmode\expandafter\endproofmath\else
  \unskip\nobreak\hfil\penalty50\hskip.75em\hbox{}\nobreak\hfil\bull
  {\parfillskip=0pt \finalhyphendemerits=0 \bigbreak}\fi}
\def\endproofmath$${\eqno\bull$$\bigbreak}
\def\bull{\vbox{\hrule\hbox{\vrule\kern3pt\vbox{\kern6pt}\kern3pt\vrule}\hrule}}
\newtheorem{theorem}{Theorem}[section]

\newtheorem{proposition}[theorem]{Proposition}
\newtheorem{lemma}[theorem]{Lemma}

\newtheorem{D}[theorem]{Definition}

\newtheorem{R}[theorem]{Remark}
\newenvironment{remark}{\begin{R}\rm }{\end{R}}

\def\Zee{\mathbb{Z}}
\def\Q{\mathbb{Q}}

\def\Cee{\mathbb{C}}

\def\Pee{\mathbb{P}}

\def\scrO{\mathcal{O}}

\def\ov{\overline}

\def\frak{\mathfrak}
\def\Pic{\operatorname{Pic}}

\def\Aut{\operatorname{Aut}}

\def\Id{\operatorname{Id}}

\def\ad{\operatorname{ad}}

\def\Lie{\operatorname{Lie}}

\title{Exceptional Groups and del Pezzo Surfaces}
\author{Robert Friedman\thanks{The first author was partially
    supported by NSF grant DMS-99-70437.}
\  and John W. Morgan\thanks{The second author was partially supported
   by NSF grant DMS-97-04507.}}

\begin{document}

\maketitle

\centerline{\it Dedicated to Herb Clemens}
\bigskip

\section*{Introduction}

Let $\mathbf{E}_r$, $r=6,7,8$ denote the simply connected form of the complex
linear group whose root system is of type $E_r$.  We extend this series to $3\leq
r\leq 8$ by setting $E_5= D_5$, $E_4=A_4$, and $E_3=A_1\times A_2$, again with the
understanding that $\mathbf{E}_r$ denotes the simply connected form of the
corresponding complex linear group of type $E_r$. It has long been known that
there are deep connections between
$\mathbf{E}_r$ and del Pezzo surfaces of degree
$9-r$. The goal of this paper is to make one of these connections explicit: if
$X$ is a del Pezzo surface of degree $d=9-r$, possibly with rational double point
singularities, we show that there is a  ``tautological"
holomorphic $\widetilde {\mathbf{E}}_r$-bundle $\Xi$ over $X$, where $\widetilde {\mathbf{E}}_r$ is an appropriate conformal form of the group $\mathbf{E}_r$. The
most classical case, $r=6$, corresponds to the case of cubic surfaces. In this
case,
$\widetilde {\mathbf{E}}_6 = \mathbf{E}_6\times _{\Zee/3\Zee}\Cee^*$. There is a
natural
$27$-dimensional representation $\rho$ of $\mathbf{E}_6$ and $\widetilde {\mathbf{E}}_6$. If
$X$ is a smooth cubic surface,  then the induced holomorphic vector bundle $\Xi 
\times_{\widetilde {\mathbf{E}}_6}\Cee^{27}$ is isomorphic to $\bigoplus
_{i=1}^{27}\scrO_X(L_i)$, where the $L_i$ are the distinct lines on $X$. The fact
that $\Xi 
\times_{\widetilde {\mathbf{E}}_6}\Cee^{27}$ is isomorphic to a direct sum of line
bundles  reflects the fact that the structure group of $\Xi$ reduces to a maximal
torus of
$\widetilde {\mathbf{E}}_6$. When 
$X$ has rational double points,  the
induced rank $27$ vector bundle is no longer a direct sum of line bundles.
Instead, the line bundle factors on a general surface coalesce into
irreducible summands of higher rank, reflecting the way in which lines coalesce
on singular cubic surfaces. Correspondingly, the structure group of $\Xi$
reduces to a reductive subgroup of $\widetilde {\mathbf{E}}_6$ whose Lie algebra 
is generated by a maximal torus and by  roots corresponding to smooth rational
curves of self-intersection $-2$ in the minimal resolution of $X$. This phenomenon
reflects the picture in physics, where rational double point singularities
correspond to extra massless particles, and these particles are described as
gauge particles for a gauge group formed in exactly the same way. Similar results
hold for any value of
$r$.

One motivation for describing the bundle $\Xi$ is to give a more direct
explanation for the correspondence described in \cite{FMW1} between S-equivalence
classes of semistable
$\mathbf{E}_r$-bundles over a smooth elliptic curve $E$ with origin $p_0$ and
triples $(X, D,
\varphi)$, where $X$ is a del Pezzo surface of degree $9-r$, $D$ is a hyperplane
section of $X$, and $\varphi\colon D \to E$ is an isomorphism from $D$ to $E$
such that $\varphi^*\scrO_E((9-r)p_0) \cong \scrO_X(D)|D$. We show that, given
the triple  $(X, D, \varphi)$, after a suitable twist by a line bundle, there is a
canonical reduction of the structure group of the bundle  $\varphi^*(\Xi|D)$ to an
$\mathbf{E}_r$-bundle $\xi$ over
$E$ which realizes the above correspondence. Moreover, in our construction, the
bundle
$\xi$ is always the regular representative.

In \cite{Manin}, Manin writes of the $27$ lines on a cubic surface, ``Their
elegant symmetry both enthrals and at the same time irritates; what use is it to
know, for instance, the number of coplanar triples of such lines (forty five) or
the number of double Schl\"affli sixfolds (thirty six)?"  One answer, different
from that given in \cite{Manin}, is that the numerology of the lines, conics,
twisted cubics, etc., is deeply related to the weights of the fundamental
representations of
$\mathbf{E}_r$. We give a brief discussion of this in Section 5.

Much of the material in this paper on del Pezzo surfaces is standard. General
references are
\cite{Manin} and \cite{777}. We should also mention that Leung has
independently given a different construction of the bundle $\Xi$ in case $X$ is
smooth \cite{Leung}.

Let us close this introduction with a list of some outstanding open questions:
\begin{enumerate}
\item Is there an algebraic construction of the bundles $\Xi$ along the lines of
the parabolic construction of bundles on the hyperplane sections \cite{FMII}?
\item In the case where $D$ is an irreducible but singular hyperplane section of
$X$, passing through a singular point, and $\varphi$ is an isomorphism from $D$
to a fixed Weierstrass cubic $E$, can one relate the bundle $\xi$ constructed
above to the $\mathbf{E}_r$-bundle constructed in \cite{FMIII}?
\item Can the construction be extended to certain more singular del Pezzo
surfaces, in particular the non-normal del Pezzo surfaces of Reid \cite{Reid}?
What about weighted cones over elliptic curves?
\item Is there a general method of constructing $\Xi$ for families of del Pezzo
surfaces $\pi\colon \mathcal{X} \to B$ over a base $B$? Does some variant of the
relative intermediate Jacobian of $\mathcal{X}$ play a role in describing all
such bundles?
\end{enumerate}

\section{Some results on lattices}

Let $\Lambda$ be a lattice. In practice, $\Lambda$ will be the coroot lattice of
the simply connected group $G$, or $\Lambda$ will be $[K_{\widetilde X}]^\perp
\subseteq H^2(\widetilde X; \Zee)$, where $\widetilde X$ is the minimal
resolution of a del Pezzo surface with rational double points. We consider
extensions
$$0 \to \Zee\oplus \Lambda \to \widetilde \Lambda \to \Zee/d\Zee \to 0,$$
where $\widetilde \Lambda$ is also a lattice.
Let $\kappa$ be the image of $(1,0)$, and identify $\Lambda$ with its image. We
assume that $\Lambda$ is a primitive sublattice of $\widetilde \Lambda$, and that
$\kappa$ is a primitive vector. There are then two related exact sequences
\begin{align*}
0 \to \Lambda \to &\widetilde \Lambda \to \Zee \to 0;\\
0 \to \Zee \to &\widetilde \Lambda \to \ov{\Lambda}\to 0.
\end{align*}
Choose a generator of $\Zee/d\Zee$, and lift it to an element $\mu \in
\widetilde \Lambda$. We may write $\mu = \mu_1+\mu_2$, where
$\mu_1\in \frac{1}{d}\Zee$ and $\mu_2\in \frac{1}{d}\Lambda$. It is easy
to check that our hypotheses imply that   $\mu_1$ has order $d$ as an element
of $\frac{1}{d}\Zee/\Zee$, and similarly for $\mu_2$. We can choose the
generator of $\Zee/d\Zee$ uniquely so that $\mu_1\equiv 1/d \mod \Zee$, and
then $\mu_2$ is a well-defined element of $\frac{1}{d}\Lambda/\Lambda$ and is a
complete invariant of the extension.

Suppose that there is a unimodular quadratic form on $\widetilde \Lambda$ such
that $\kappa$ and $\Lambda$ are orthogonal. In this case, $\mu_2$ is a
well-defined element of $\Lambda^*/\Lambda$. Given $\lambda \in \widetilde
\Lambda$, we can write $\lambda =\lambda_1 + \lambda_2$, where the 
$\lambda_i\in\widetilde
\Lambda\otimes \Q$, $\lambda_1 = \langle \lambda, \kappa
\rangle\kappa/\kappa^2$, and $\lambda_2\in \kappa^\perp$.  Thus $\lambda_i\in
\widetilde \Lambda$ if and only if 
$\kappa^2|\langle \lambda, \kappa
\rangle$, and the index $d$ of $\Zee\oplus \Lambda$ in $\widetilde \Lambda$
is $|\kappa^2|$. We can choose a representative for
$\mu_2$ as follows: first choose $\mu \in \widetilde \Lambda$ such that
$\langle \mu, \kappa \rangle = 1$. Then we can take
$$\mu _2 \equiv \mu - \frac{1}{\kappa^2}\kappa \mod \Lambda.$$
In fact, $\Lambda^*/\Lambda$ is cyclic of order $d$ and $\mu_2$ is a generator.

\medskip
\noindent \textbf{Example:}  Suppose that $3\leq r \leq 8$. If $\widetilde 
\Lambda =\Zee^{r+1}$ with diagonal basis
$h, e_1,
\dots, e_r$, where $h^2=1$ and $e_i^2=-1$, let $\kappa = 3h-\sum _ie_i$ and
$\Lambda =\kappa^\perp$. Take $\mu = e_r$. A basis for $\Lambda$ is given by
$\alpha_1 = e_1-e_2, \dots, \alpha_{r-1} = e_{r-1}-e_r, \alpha _r =
3h-e_1-e_2-e_3$. In this case $\mu_2 =\varpi_{r-1}$ is the dual basis element
to $\alpha_{r-1}$. 

The lattice $\Lambda$ is the coroot lattice for $\mathbf{E}_r$, with the negative
of the usual intersection form. The action of the Weyl group $W$ generated by
reflections in the roots extends to an action of $W$ on $\widetilde \Lambda$. We
have the following:

\begin{lemma}\label{autom} The group $W$ is exactly the subgroup of automorphisms
of
$\widetilde \Lambda$ which fix $\kappa$.
\end{lemma}
\begin{proof} The automorphism group of $\Lambda$ is a semidirect product of $W$
and the outer automorphism group of the root system corresponding to $\Lambda$.
It is easy to verify that, in the cases at hand, a nontrivial  outer automorphism
acts by sending $\mu_2$ to $-\mu_2 \mod \Lambda$. From this, the result is clear.
\end{proof}

The intersection form identifies $\widetilde \Lambda$ with its dual. Using this
identification, we have exact sequences
\begin{align*}
0 \to \Zee\to &\widetilde \Lambda \to \Lambda ^* \to 0;\\
0 \to \widetilde \Lambda \to (\Zee&\oplus \Lambda)^* \to \Zee/d\Zee\to 0.
\end{align*}
From these sequences, we see the following:

\begin{lemma}\label{duals} 
\begin{enumerate}
\item[\rm (i)] A homomorphism $\phi\colon \Zee\oplus \Lambda \to \Zee$ is given
by inner product with an element of $\widetilde \Lambda$ if and only if
$\phi(1,0) \equiv -\phi(0,d\mu_2) \mod d$.
\item[\rm (ii)] Given $\phi\in \Lambda ^*$, there exists a $\lambda \in \Lambda$
such that $\phi(v) = \langle\lambda, v\rangle$ for all $v\in \Lambda$. Moreover,
any two choices of $\lambda$ differ by an integral multiple of $\kappa$. \qed
\end{enumerate}
\end{lemma}

\section{A tautological bundle over the Cartan subgroup of a conformal form}

We begin by constructing the appropriate conformal form of the group
$\mathbf{E}_r$. More generally, suppose that $G$ is a reductive group with Cartan
subgroup $H$, and that $c$ is an element of the center of $G$ of finite order
$d$. Let $\Lambda =
\pi_1(H)$. We can then form the group $\widetilde G = \Cee^*\times
_{\Zee/d\Zee}G$, where the generator
$1\in
\Zee/d\Zee$ maps to $\exp(2\pi\sqrt{-1}/d) \in \Cee^*$ and to $c\in G$. A Cartan
subgroup of $\widetilde G$ is given by $\widetilde H = \Cee^*\times
_{\Zee/d\Zee}H$. Let
$\widetilde \Lambda =\pi_1(\widetilde H)$. 
From the exact sequence
$$0 \to \Zee/d\Zee \to \Cee^* \times H \to \Cee^*\times _{\Zee/d\Zee}H \to 0,$$
we see that there is an exact sequence
$$0\to \Zee \oplus \Lambda \to \widetilde \Lambda \to \Zee/d\Zee \to 0,$$
and that we obtain the extension of tori by tensoring the above sequence
with $\Cee^*$. It is an easy exercise to check that the element $\mu_2$ defined
in the previous section is just the element $c^{-1}$, viewed as an element of
$\Lambda^*/\Lambda\subseteq \Lambda\otimes \Q/\Lambda$. Conversely, given an
extension $\widetilde \Lambda$ of $\Zee \oplus \Lambda$ with a unimodular form as
in the previous section, suppose that the element $\mu_2$ corresponds to the
element $c^{-1}\in \Lambda^*/\Lambda$. If we tensor with
$\Cee^*$, we get an extension of tori and the torus $\widetilde H=\widetilde
\Lambda \otimes _\Zee\Cee^*$ is a Cartan
subgroup of  $\widetilde G$.

We apply this to the lattice $\widetilde\Lambda =\Zee^{r+1}$ with basis $h, e_1,
\dots, e_r$ of the example given at the end of the previous section. In this case,
$\Lambda$ is identified with the coroot lattice of the group $\mathbf{E}_r$ (with
the negative of the usual intersection form), and thus $\widetilde H=\widetilde
\Lambda \otimes _\Zee\Cee^*$ is a Cartan subgroup for
$\widetilde {\mathbf{E}}_r = \mathbf{E}_r\times _{\Zee/(9-r)\Zee}\Cee^*$. Note
that, for
$3\leq r\leq 8$ and $r\neq 4$, the group $\Zee/(9-r)\Zee$ has at most two
generators and there is essentially no choice for the group $\widetilde {\mathbf{E}}_r$. However, for $r=4$, the central element $c$ is the square of the
standard generator of the center of
$SL_5(\Cee)$, and $\widetilde {\mathbf{E}}_4$ is not $GL_5(\Cee)$ but rather the
double cover of $GL_5(\Cee)$.

Now let  $X$ be  a del Pezzo surface of
degree $9-r$ with at worst rational double point singularities and let
$\widetilde X$ be the minimal resolution of $X$. Fix an isomorphism
$\psi\colon H^2(\widetilde X; \Zee) \to \widetilde \Lambda$ such that
$\psi ([K_X]) =-\kappa$. We will refer to such a $\psi$ as a \textsl{marking}. It
follows from Lemma~\ref{autom} that any two markings differ by an action of the
Weyl group on
$\widetilde \Lambda$. 

Holomorphic $\widetilde H$-bundles over $\widetilde X$ are classified by
$H^1(\widetilde X; \underline{\widetilde H})$, where $\underline{\widetilde H}$
is the sheaf whose sections over an open set $U$ are holomorphic maps from  $U$
to $\widetilde H$. From the exact sequence
$$0 \to \widetilde \Lambda \to \Lie (\widetilde H)\otimes \scrO_{\widetilde X}
\to \underline{\widetilde H} \to \{1\},$$ and the fact that $H^i(\widetilde X;
\scrO_{\widetilde X}) = 0$ for $i>0$, we see that there are canonical
isomorphisms 
$$H^1(\widetilde X;\underline{\widetilde H}) \cong H^2(\widetilde X;\widetilde
\Lambda) \cong H^2(\widetilde X;\Zee) \otimes _\Zee\widetilde \Lambda.$$
The marking $\psi$ defines an isomorphism $H^2(\widetilde X;\Zee) \otimes
_\Zee\widetilde \Lambda \cong \widetilde \Lambda\otimes
_\Zee\widetilde \Lambda$. Using the intersection form to identify $\widetilde
\Lambda$ with its dual, we can take the class of $\Id \in  \widetilde
\Lambda^*\otimes _\Zee\widetilde \Lambda \cong \widetilde \Lambda\otimes
_\Zee\widetilde \Lambda$. Let $\Xi_{\widetilde H}$ be the corresponding
$\widetilde H$-bundle over
$\widetilde X$, and   define
$$\Xi_{\textrm{split}}= \Xi_{\widetilde H}\times _{\widetilde H}\widetilde 
{\mathbf{E}}_r.$$ The bundle  $\Xi_{\textrm{split}}$ is split, in the sense that,
given any representation $\rho\colon \widetilde {\mathbf{E}}_r\to GL(N)$, the
induced vector bundle 
$\Xi_{\textrm{split}}\times _{\widetilde {\mathbf{E}}_r}\Cee^N$ is a direct sum of
line bundles. Since any two markings differ by an element of the Weyl group, it is
clear that the isomorphism class of the
$\widetilde {\mathbf{E}}_r$-bundle
$\Xi_{\textrm{split}}$ is independent of the choice of a marking. In case
$\widetilde X = X$, i.e. in the case of a smooth del Pezzo surface, we take
$\Xi =
\Xi_{\textrm{split}}$. In the next section, we will discuss how to modify
$\Xi_{\textrm{split}}$ so as to make it trivial in a neighborhood of the
preimages of the double points on $X$.

Let us record one basic property of $\Xi_{\textrm{split}}$. Suppose that $\lambda
\in H^2(\widetilde X;\Zee)$. Using the marking to view $\psi$ as an
element of $\widetilde
\Lambda$ and the quadratic form, we can view
$\lambda$ as a homomorphism $\widetilde \Lambda \to \Zee$. Hence $\lambda$
defines a homomorphism $\widetilde H\to \Cee^*$. There is the associated
$\Cee^*$-bundle $\Xi_{\textrm{split}}\times _{\widetilde H}\Cee^*$, which
corresponds to a complex line bundle $L_\lambda$ on $\widetilde X$. Unwinding the
definitions gives:

\begin{lemma}\label{taut} In the above notation, $c_1(L_\lambda)
=\lambda$.\qed
\end{lemma}

\section{Extension over the Borel subgroup}

The bundle $\Xi_{\textrm{split}}$ is the correct bundle only when $X=\widetilde
X$, i.e. $X$ is smooth. In case $X$ is singular, there are rational curves $C$ on
$\widetilde X$ of self-intersection
$-2$, which for brevity we call \textsl{$-2$-curves}. Define the
\textsl{exceptional divisor} $Y$ on $\widetilde X$ to be the union of the
$-2$-curves. The bundle $\Xi_{\textrm{split}}$ is nontrivial in every
neighborhood of a $-2$-curve, and hence does not descend to a bundle on $X$.   
In this section, we
show how to modify the bundle
$\Xi_{\textrm{split}}$ so that it becomes trivial in a neighborhood of the
exceptional divisor. The idea is as follows. Fix
a general ample line bundle $L$ on $\widetilde X$. Then $\psi(L)$ defines a
positive Weyl chamber in $\Lambda$ and hence a Borel subgroup of $\mathbf{E}_r$
and thus of $\widetilde {\mathbf{E}}_r$. Let $U$ be the unipotent radical of this
Borel subgroup, so that the Borel is equal to $\widetilde H U$. We will construct
a lift of the
$\widetilde H$-bundle
$\Xi_{\widetilde H}$ to an
$\widetilde H U$-bundle $\Xi_{\widetilde H U}$, so that $\Xi_{\widetilde H
U}\times _{\widetilde H U}\widetilde {\mathbf{E}}_r$ is trivial in an analytic 
neighborhood of the exceptional divisor, and hence arises from a bundle over
$X$.  The results are summarized by the following theorem:

\begin{theorem}\label{extend} There exists a unique $\widetilde {\mathbf{E}}_r$-bundle
$\widetilde
\Xi$ over
$\widetilde X$ with the following properties:
\begin{enumerate}
\item[\rm (i)] There exists a bundle $\Xi_{\widetilde H U}$, where $U$ is the
unipotent radical of the Borel subgroup of $\widetilde {\mathbf{E}}_r$, such that
$\Xi_{\widetilde H U}/U =\Xi_{\rm split}$ and $\widetilde \Xi = \Xi _{\widetilde
H U} \times _{\widetilde H U}\widetilde {\mathbf{E}}_r$.
\item[\rm (ii)] The bundle $\widetilde \Xi$ is rigid, i.e. $H^1(\widetilde X; \ad
\widetilde \Xi) =0$.
\end{enumerate}
In case $X=\widetilde X$, $\widetilde \Xi = \Xi_{\rm split}$. Furthermore, the
bundle $\widetilde \Xi$ has the following properties:
\begin{enumerate}
\item[\rm (iii)] For every $D\in |-K_{\widetilde X}|$, $h^0(D;  \ad
\widetilde \Xi|D) = r+1$.
\item[\rm (iv)] There exists a unique $\widetilde {\mathbf{E}}_r$-bundle $\Xi$ on
$X$ such that $\widetilde \Xi =p^*\Xi$. 
\end{enumerate}
\end{theorem}

The proof will be in several steps. We begin with some very general and well-known
results about
$-2$-curves on the minimal resolution of a del Pezzo surface $X$ with rational
double points. Let
$p\colon \widetilde X\to X$ be the minimal resolution of $X$. Given $\alpha \in
H^2(\widetilde X;\Zee)$, we call $\alpha$ a \textsl{root} if $\alpha^2=-2$ and
$\alpha \cdot K_{\widetilde X} = 0$. For any class $\alpha \in H^2(\widetilde
X;\Zee)$, let
$L_\alpha$ be the corresponding line bundle.

\begin{lemma}\label{31} Let $\alpha$ be a root,  let $C_1, \dots, C_k$ be $-2$
curves in $\widetilde X$, and let $a_i\in \Zee$. Then $$H^2(\widetilde X;
L_\alpha\otimes
\scrO_{\widetilde X}(\sum _ia_iC_i)) = 0.$$
 Moreover,  the following are
equivalent:
\begin{enumerate}
\item[\rm (i)] $H^1(\widetilde X; L_\alpha)\neq 0$;
\item[\rm (ii)] $H^1(\widetilde X; L_\alpha)\cong \Cee$;
\item[\rm (iii)] $H^0(\widetilde X; L_\alpha)\neq 0$, i.e. $L_\alpha =
\scrO_{\widetilde X}(C)$ for some effective curve
$C$;
\item[\rm (iv)] $H^0(\widetilde X; L_\alpha)\cong \Cee$.
\end{enumerate}
\end{lemma}
\begin{proof} By Serre duality, $H^2(\widetilde X; L_\alpha)\cong H^0(\widetilde
X; K_{\widetilde X}\otimes L_\alpha^{-1})$. Since $K_{\widetilde X}
=\scrO_{\widetilde X}(-D)$, where $D$ is nef and big, and $D\cdot \alpha = 0$,
$H^2(\widetilde X; L_\alpha) =0$ for every root $\alpha$. A similar statement
holds for $H^2(\widetilde X; L_\alpha\otimes \scrO_{\widetilde
X}(\sum _ia_iC_i))$.

Let $\alpha$ be a root. Then $\alpha ^2=-2$, and so by the Riemann-Roch theorem,
$\chi(\widetilde X; L_\alpha) = 0$. By the first part of the lemma,
$h^2(\widetilde X; L_\alpha) = 0$ and so
$h^1(\widetilde X; L_\alpha) = h^0(\widetilde X; L_\alpha)$. Thus $h^1(\widetilde
X; L_\alpha) \neq 0$ if and only if $L_\alpha = \scrO_{\widetilde X}(C)$ for some 
effective curve
$C$. In this case, since $D\cdot C =0$ and $D$ is nef and big, it follows that
$h^0(\widetilde X; L_\alpha) = 1$, and hence that $h^1(\widetilde X; L_\alpha)
=1$. The remaining equivalences are clear.
\end{proof}

If $\alpha$ satisfies any of the equivalent conditions above, we call $\alpha$ an
\textsl{effective  root}.

Next, we see how the line bundles $L_\alpha$ restrict to smooth elements $D\in
|-K_{\widetilde X}|$:

\begin{lemma}\label{restricttoD} Let $\alpha$ be a root. Suppose that $D$ is
reduced and irreducible and that
$K_{\widetilde X} =\scrO_{\widetilde X}(-D)$. Then the following are equivalent:
\begin{enumerate}
\item[\rm (i)] $L_\alpha|D =\scrO_D$;
\item[\rm (ii)] Either $\alpha$ or $-\alpha$ is effective.
\end{enumerate}
Moreover, in this case $\alpha$ is effective if and only if $L\cdot \alpha > 0$.
Finally, for every root $\alpha$, if
$L\cdot \alpha > 0$, then the natural maps
$H^i(\widetilde X; L_\alpha) \to H^i(D; L_\alpha|D)$ are all isomorphisms.
\end{lemma}
\begin{proof} In any case $L_\alpha|D$ has degree zero. Suppose that $-\alpha$ is
not effective. Then
$H^i(\widetilde X; L_\alpha^{-1}) = 0$ for all $i$. Consider the long exact
cohomology sequence arising from
$$0\to \scrO_{\widetilde X}(-D) \otimes L_\alpha \to L_\alpha\to L_\alpha |D \to
0.$$ Since $H^i(\widetilde X;\scrO_{\widetilde X}(-D) \otimes L_\alpha)$ is Serre
dual to $H^{2-i}(\widetilde X;L_\alpha^{-1})$, it is zero for all $i$. Thus 
 $H^i(\widetilde X;
L_\alpha) \to H^i(D; L_\alpha|D)$ is an isomorphism for all $i$. It follows that
$\alpha$ is effective if and only if $H^0(D; L_\alpha|D)\neq 0$, if and only if
$L_\alpha|D$ is trivial (since it is a line bundle of degree zero). The remaining
statements are clear.
\end{proof}

The next lemma gives a necessary and sufficient condition for the bundle
$\widetilde \Xi$ to descend to a bundle on the singular surface $X$.

\begin{lemma}\label{descend} Let $\widetilde \Xi$ be an
$\widetilde{\mathbf{E}}_r$-bundle over
$\widetilde X$. Suppose that $\widetilde \Xi|C$ is the trivial $\widetilde
{\mathbf{E}}_r$-bundle for every irreducible $-2$-curve $C$.
Then there exists a unique $\widetilde {\mathbf{E}}_r$-bundle $\Xi$ on $X$ such that
$\widetilde
\Xi = p^*\Xi$.
\end{lemma}
\begin{proof} Let $V$ be  a vector bundle on $\widetilde X$  such that $V|C \cong
\scrO_C^n$ for every irreducible $-2$-curve $C$. We shall show that $V$ is
trivial in some analytic neighborhood of the exceptional curve $Y$. If
$C_1,
\dots, C_e$ are   the
$-2$-curves, and $Y = \bigcup_iC_i$, it follows easily from the fact that the
dual graph of
$Y$ is a union of contractible components and that the $C_i$ meet transversally
that $V|Y =\scrO_Y^n$. Let
$U$ be a contractible Stein neighborhood of
$X_{\textrm{sing}}$, possibly disconnected, and let $\widetilde U = p^{-1}(U)$. 
We claim that the natural map $H^0(\widetilde U; V|\widetilde U)\to H^0(Y; V|Y)$
is surjective. In fact, the cokernel of this map is contained in $H^1(\widetilde
U; V\otimes \scrO_{\widetilde X}(-Y)|\widetilde U) \cong H^0(R^1 p_*(V\otimes
\scrO_{\widetilde X}(-Y)))$. A standard argument using the fact that the
singularities are rational double points and the formal functions theorem shows
that $R^1 p_*(V\otimes
\scrO_{\widetilde X}(-Y)) =0$. Thus $H^0(\widetilde U; V|\widetilde U)\to H^0(Y;
V|Y)$ is surjective. Lifting a basis of sections of $V|Y$ to
$\widetilde U$, we can assume, possibly after shrinking $U$ and $\widetilde U$,
that $V|\widetilde U$ is trivial.

Now choose a faithful representation $\rho$ of $\widetilde
{\mathbf{E}}_r$, which we use to view  
$\widetilde {\mathbf{E}}_r$ as a subgroup of $GL(N)$ for some  $N$. Let $V$ be 
the induced vector bundle $\widetilde
\Xi\times _{\widetilde {\mathbf{E}}_r}\Cee^N$.  By assumption, $V|C$ is the
trivial vector bundle
$\scrO_C^N$ for every $-2$-curve $C$. Thus, $V|\widetilde U$ is the trivial
vector bundle over $\widetilde U$.

Let $Q$ be
the quotient space
$GL(N)/\widetilde {\mathbf{E}}_r$, which exists as an affine variety (since
$\widetilde {\mathbf{E}}_r$ is reductive) and hence as a complex manifold. Let
$\underline{Q}$ be the sheaf of morphisms from $\widetilde {U}$ to $Q$, and
similarly for
$\underline{\widetilde {\mathbf{E}}_r}$ and $\underline{GL(N)}$. From the exact
sequence of pointed sets
$$\{1\} \to \underline{\widetilde {\mathbf{E}}_r} \to \underline{GL(N)} \to
\underline{Q}\to
\{1\},$$ we have a long exact cohomology sequence through the term
$H^1(\widetilde{U}; \underline{GL(N)})$. If
$\ov\xi$ is the class in $H^1(\widetilde U; \underline{\widetilde
{\mathbf{E}}_r} )$ induced by
$\widetilde \Xi|\widetilde U$, then $\ov \xi$ maps to the trivial element in
$H^1(\widetilde U; 
\underline{GL(N)})$, and hence comes by coboundary from $H^0(\widetilde U;
\underline{Q})$, which we can identify with the set of morphisms from
$\widetilde U$ to $Q$. If
$f$ is such a morphism, then since $Q$ is affine, $f(Y) =\textrm{pt}$. Since $U$
is normal, it follows that $f$ is induced from a morphism $U\to Q$. But since
there are local cross sections for the map $GL(N)
\to Q$ (in the classical topology), again after shrinking $U$ we can assume that
$f$ is in the image of a morphism from $\widetilde U$ to $GL(N)$. In particular,
this says that $\ov \xi$ is trivial. Thus the triviality of $V|\widetilde U$
implies that of
$\widetilde
\Xi|\widetilde U$. Choosing a local trivialization of $\widetilde \Xi$ for which
$\widetilde U$ is one of the open sets and the remaining ones do not meet $Y$
gives a $1$-cocycle defining $\widetilde\Xi$ which also defines a bundle $\Xi$
on $X$, and clearly $\widetilde \Xi =
p^*\Xi$. The uniqueness is a straightforward consequence of the fact that, as
$\widetilde {\mathbf{E}}_r$ is affine, every morphism from $\widetilde U$ into
$\widetilde {\mathbf{E}}_r$ is constant on the exceptional curves and hence
descends to a morphism from
$U$ into $\widetilde {\mathbf{E}}_r$.
\end{proof}

We turn now to the construction of a suitable bundle $\widetilde \Xi$, which we
shall always take to be of the form $\Xi_{\widetilde H
U}\times _{\widetilde H U}\widetilde {\mathbf{E}}_r$, where $\Xi_{\widetilde H U}$
is a lift of the $\widetilde H$-bundle
$\Xi_{\widetilde H}$ to an
$\widetilde H U$-bundle.  

\begin{lemma}\label{rigid} Let $\widetilde \Xi = \Xi_{\widetilde H
U}\times _{\widetilde H U}\widetilde {\mathbf{E}}_r$, where $\Xi_{\widetilde H U}$
is a lift of the $\widetilde H$-bundle
$\Xi_{\widetilde H}$ to an
$\widetilde H U$-bundle. Then the following are equivalent:
\begin{enumerate}
\item[\rm (i)] $H^1(\widetilde X; \ad \widetilde \Xi) = 0$;
\item[\rm (ii)] $h^0(\widetilde X; \ad \widetilde \Xi) = r+1$;
\item[\rm (iii)] There exists a smooth  $D\in |-K_{\widetilde X}|$ such that 
$h^0(D; \ad \widetilde \Xi|D) = r+1$;
\item[\rm (iv)] For every  $D\in |-K_{\widetilde X}|$,
$h^0(D; \ad \widetilde \Xi|D) = r+1$.
\end{enumerate}
\end{lemma}
\begin{proof} First we claim:

\begin{lemma}\label{h2zero} With $\widetilde \Xi$ as in the statement of
Lemma~\ref{rigid},
\begin{enumerate} 
\item[\rm (i)] $\chi(\widetilde X; \ad \widetilde \Xi)= r+1$;
\item[\rm (ii)] Let $C_1, \dots, C_k$ be $-2$-curves and let $a_i\in \Zee$. Then
$$H^2(\widetilde X; \ad \widetilde \Xi\otimes \scrO_{\widetilde X}(\sum _ia_iC_i))
= 0 = H^0(\widetilde X; \ad \widetilde \Xi\otimes \scrO_{\widetilde X}(\sum
_ia_iC_i)\otimes \scrO_{\widetilde X}(-D)) = 0.$$
\end{enumerate}
\end{lemma}
\begin{proof} (i) By Riemann-Roch, $\chi(\widetilde X; \ad \widetilde \Xi)$ only
depends on the topological type of $\widetilde \Xi$. Thus we may replace
$\widetilde \Xi$ by $\Xi_{\textrm{split}}$. But
$$\ad\Xi_{\textrm{split}} =\scrO_{\widetilde X}^{r+1} \oplus \bigoplus _{\alpha
\in R}L_\alpha,$$
and as we have seen, $\chi(\widetilde X; L_\alpha)=0$ for every root $\alpha$.
Hence $\chi(\widetilde X; \ad \widetilde \Xi)= (r+1)\chi(\widetilde
X;\scrO_{\widetilde X}) = r+1$.

(ii) The bundle $\ad \widetilde \Xi$ has a filtration whose successive quotients
are either $\scrO_{\widetilde X}$ or $L_\alpha$, $\alpha \in R$. Thus, by
Lemma~\ref{31}, $H^2(\widetilde X; \ad \widetilde \Xi\otimes  \scrO_{\widetilde
X}(\sum _ia_iC_i)) = 0$, and a similar argument shows that $H^0(\widetilde X; 
\ad \widetilde \Xi\otimes  \scrO_{\widetilde X}(\sum _ia_iC_i)\otimes
\scrO_{\widetilde X}(-D)) = 0$.
\end{proof}

Returning to the proof of Lemma~\ref{rigid}, we see that 
$$h^0(\widetilde X; \ad
\widetilde \Xi) - h^1(\widetilde X; \ad \widetilde \Xi) =\chi(\widetilde X; \ad
\widetilde \Xi) =r+1.$$
Thus $h^0(\widetilde X; \ad \widetilde \Xi) \geq r+1$,
with equality if and only if $h^1(\widetilde X; \ad \widetilde \Xi)=0$. In
particular, (i) $\iff$ (ii). Let $D\in |-K_{\widetilde X}|$. From the exact
sequence
$$0 \to \ad \widetilde \Xi \otimes \scrO_{\widetilde X}(-D) \to \ad \widetilde
\Xi \to \ad \widetilde \Xi|D \to 0,$$
and the fact that $H^1(\widetilde X; \ad\widetilde \Xi\otimes \scrO_{\widetilde
X}(-D))$ is Serre dual to $H^1(\widetilde X; \ad\widetilde \Xi)$, it follows that,
if $h^1(\widetilde X; \ad \widetilde \Xi)=0$, then $h^0(D;\ad \widetilde \Xi|D) =
 h^0(\widetilde X; \ad
\widetilde \Xi)=r+1$. Thus (i) $\implies$ (iv)
$\implies$ (iii). Finally, since $H^0(\widetilde X; \ad\widetilde \Xi\otimes
\scrO_{\widetilde X}(-D))=0$, the map $H^0(\widetilde X; \ad\widetilde \Xi) \to
H^0(D; \ad\widetilde \Xi|D)$ is injective. As we have seen, $h^0(\widetilde X;
\ad\widetilde \Xi) \geq r+1$. Thus, if $h^0(D; \ad\widetilde \Xi|D) =r+1$, then 
$h^0(\widetilde X;
\ad\widetilde \Xi) = r+1$ as well, so that (iii) $\implies$ (ii).
\end{proof}

Now let us analyze the restriction of a bundle $\widetilde \Xi$ satisfying any of
the equivalent hypotheses of Lemma~\ref{rigid} to a $-2$-curve.

\begin{lemma}\label{trivial} Let $\widetilde \Xi = \Xi_{\widetilde H
U}\times _{\widetilde H U}\widetilde {\mathbf{E}}_r$, where $\Xi_{\widetilde H U}$
is a lift of the $\widetilde H$-bundle
$\Xi_{\widetilde H}$ to an
$\widetilde H U$-bundle. Suppose that $H^1(\widetilde X; \ad \widetilde \Xi)=0$.
Then, for every $-2$-curve $C$, $\widetilde \Xi|C$ is the trivial  $\widetilde
{\mathbf{E}}_r$-bundle. 
\end{lemma}
\begin{proof} First, we claim that $\widetilde \Xi|C$ is topologically trivial.
It suffices to check that $\Xi_{\textrm{split}}|C$ is topologically trivial. The
surjection $\widetilde {\mathbf{E}}_r\to \Cee^*$ induces a $\Cee^*$-bundle $\det
(\Xi_{\textrm{split}}|C)$, and the corresponding line bundle is easily seen to be
$K_{\widetilde X}|C$, which is trivial. Hence $\Xi_{\textrm{split}}|C$ lifts to a
$\mathbf{E}_r$-bundle, which is automatically topologically trivial since
$\mathbf{E}_r$ is simply connected.

We have the exact sequence 
$$0 = H^1(\widetilde X; \ad \widetilde \Xi) \to H^1(C; \ad\widetilde \Xi|C)
 \to H^2(\widetilde X; \ad \widetilde \Xi\otimes \scrO_{\widetilde X}(-C)).$$
By Lemma~\ref{h2zero}, $H^2(\widetilde X; \ad \widetilde \Xi\otimes
\scrO_{\widetilde X}(-C))=0$. Hence $H^1(C; \ad\widetilde \Xi|C) =0$ as well, so
that $\widetilde \Xi|C$ is rigid. But, since the trivial bundle is semistable,
the only rigid $\widetilde {\mathbf{E}}_r$-bundle on $C$ is the trivial bundle. Thus
$\widetilde
\Xi|C$ is trivial.
\end{proof}

Next we  construct the bundle $\widetilde \Xi$:

\begin{lemma}\label{overBorel}   There
exists an
$\widetilde H U$-bundle $\Xi_{\widetilde H U}$ lifting $\Xi_{\widetilde
H}$ such that, for $\widetilde \Xi =
\Xi_{\widetilde H U}\times _{\widetilde H U}\widetilde {\mathbf{E}}_r$, and for
every 
$D\in |-K_{\widetilde X}|$, we have 
$h^0(D;\ad
\widetilde \Xi|D) = r+1$. 
\end{lemma}
\begin{proof} By Lemma~\ref{rigid}, it suffices to prove the result for one
choice of $D\in |-K_{\widetilde X}|$, which we may assume to be a smooth elliptic
curve.
 By \cite{FMI}, there exists a
$\widetilde {\mathbf{E}}_r$-bundle $\widetilde \xi$ of the form $\xi_{\widetilde H
U}\times _{\widetilde H U}\widetilde {\mathbf{E}}_r$, such that (i) $h^0(D; \ad
\widetilde
\xi) = r+1$, and (ii)
$\xi_{\widetilde H U}$ is a lift to $\widetilde H U$ of the bundle
$\xi_{\widetilde H} = \Xi_{\widetilde H}|D$. Let us show that we can find an
$\widetilde H U$-bundle $\Xi_{\widetilde H U}$ such that $\Xi_{\widetilde H U}|D
= \xi_{\widetilde H U}$.  To make the construction, choose any decreasing 
filtration $U=U_0\supseteq U_1  \supseteq \cdots \supseteq U_N\supseteq
U_{N+1}=\{1\}$  of
$U$ by normal
$\widetilde H$-invariant subgroups such that, for every $i$, 
$U_i/U_{i+1}$ is contained in the center of $U/U_{i+1}$. It follows that
$U_i/U_{i+1}$ is a vector group which is a direct sum of root spaces $\frak
g^\alpha$, $\alpha \in R_i$, say, and $\widetilde H$ acts on $U_i/U_{i+1}$ in the
usual way. By  induction on
$i$, starting with
$i=0$, it suffices to show the following: Given
$i$, suppose that we have found an $\widetilde HU/U_i$-bundle $\Xi_{\widetilde
HU/U_i}$ such that   
$\Xi_{\widetilde HU/U_i}|D\cong \xi_{\widetilde H U}/U_i$. Then we can
lift
$\Xi_{\widetilde HU/U_i}$ to a $\widetilde HU/U_{i+1}$-bundle $\Xi_{\widetilde
HU/U_{i+1}}$ such that $\Xi_{\widetilde HU/U_{i+1}}|D\cong \xi_{\widetilde H
U}/U_{i+1}$. By general formalism (we use the notation of \cite[Appendix]{FMII}),
the obstruction to lifting
$\Xi_{\widetilde HU/U_i}$ to some $\widetilde HU/U_{i+1}$-bundle lives in the
cohomology group $H^2(\widetilde X; (U_i/U_{i+1})(\Xi_{\widetilde HU/U_i}))$. But
$$(U_i/U_{i+1})(\Xi_{\widetilde HU/U_i}))\cong \bigoplus _{\alpha \in
R_i}L_\alpha,$$
and $H^2(\widetilde X;L_\alpha) = 0$. Thus there is a lift of  $\Xi_{\widetilde
HU/U_i}$ to some $\widetilde HU/U_{i+1}$-bundle $\Xi_{\widetilde HU/U_{i+1}}$.
Let $\xi_{\widetilde HU/U_{i+1}} = \Xi_{\widetilde HU/U_{i+1}}|D$. Then
$\xi_{\widetilde HU/U_{i+1}}$ and $\xi_{\widetilde H
U}/U_{i+1}$ are two lifts of $\xi_{\widetilde H U}/U_i$ to $\widetilde
HU/U_{i+1}$, and as such they differ by the action of
$H^1(D;(U_i/U_{i+1})(\xi_{\widetilde H U}/U_i)$. On the other hand,
$$(U_i/U_{i+1})(\xi_{\widetilde H U}/U_i) \cong \bigoplus _{\alpha \in
R_i}(L_\alpha |D).$$
By the last sentence of Lemma~\ref{restricttoD}, since all of the roots in $R_i$
are positive, the restriction map
$H^1(\widetilde X;L_\alpha) \to H^1(D; L_\alpha |D)$ is an isomorphism for every
$\alpha \in R_i$. Thus, we can adjust the lift $\Xi_{\widetilde HU/U_{i+1}}$ by
the action of $H^1(\widetilde X; (U_i/U_{i+1})(\Xi_{\widetilde HU/U_i}))$ so that
$\Xi_{\widetilde HU/U_{i+1}}|D\cong \xi_{\widetilde H
U}/U_{i+1}$. This completes the inductive step. 
\end{proof}

Parts (i)--(iii) of Theorem~\ref{extend} now follow from Lemma~\ref{overBorel}
and Lemma~\ref{rigid}.  Part (iv) follows from Lemma~\ref{trivial}  and
Lemma~\ref{descend}.  Finally, we must check the uniqueness statement, that any
bundle satisfying (i) and (ii) of Theorem~\ref{extend} is isomorphic to
$\widetilde \Xi$. Since
$H^2(\widetilde X; \ad \Xi_{\textrm{split}}) = 0$, the deformations of
$\Xi_{\textrm{split}}$ are unobstructed, and there exists a germ of a complex
manifold, $\mathcal{U}$ which is the base space of a local semiuniversal
deformation of $\Xi_{\textrm{split}}$. By a standard argument (cf.\
\cite[Lemma 4.1.1]{FMII} for a related construction), given any extension of
$\Xi_{\textrm{split}}$ to an $\widetilde H U$-bundle $\widehat \Xi$, there is a
family of bundles over $\widetilde X\times \Cee$ which restrict to $\widehat
\Xi\times _{\widetilde H U}\widetilde{\mathbf{E}}_r$ over $\Cee-\{0\}$ and to
$\Xi_{\textrm{split}}$ over
$0$, and hence there are points of $\mathcal{U}$ corresponding to $\widehat
\Xi\times _{\widetilde H U}\widetilde{\mathbf{E}}_r$.
The subset
$\mathcal{U}'$ of
$\mathcal{U}$ defined by bundles $\Upsilon$ such that $H^1(\widetilde X; \ad
\Upsilon) = 0$ is the complement of a proper analytic subvariety and hence is
open and dense and therefore connected. Since $\widetilde \Xi$ is rigid, the
subset of $\mathcal{U}'$ corresponding to bundles isomorphic to $\widetilde \Xi$
is open, as is its complement in $\mathcal{U}'$. Hence every point of
$\mathcal{U}'$ corresponds to $\widetilde \Xi$, proving the uniqueness.
\qed

\begin{remark} The proof actually shows the following. Let $G$ be the reductive
subgroup of $\widetilde {\mathbf{E}}_r$ whose maximal torus is $\widetilde H$ and
whose Lie algebra contains all of the root spaces corresponding to the effective
roots. Then
$\widetilde \Xi$ reduces to a Borel subgroup $\widetilde H U'$ of $G$, where $U'$
is the unipotent radical of the Borel, and $\Xi$ reduces to a $G$-bundle.
\end{remark}

\section{Restriction to hyperplane sections}

We now relate the restriction of this construction to smooth hyperplane sections
to the correspondence outlined in \cite[\S2]{FMW1}. Let $X$ be a del Pezzo
surface  of degree $r$ with at worst rational double points and let $D$ be a
smooth hyperplane section. Let
$\varphi\colon D
\to E$ be an isomorphism such that $\varphi^*\scrO_E((9-r)p_0) = \scrO_X(D)|D$.
Then, following \cite[\S2]{FMW1}, we define the \textsl{period point} of $(X,D,
\varphi)$  in
$(E\otimes_\Zee \Lambda)/W$ as follows: Fix a marking $\psi$ and use it to
identify $\widetilde \Lambda$ with $\Pic \widetilde X$. Given
$\alpha\in \Lambda$, the line bundle $(\varphi^*)^{-1}(L_\alpha)$ is a line
bundle of degree zero on
$E$, and hence defines a point in $\Pic^0E\cong E$. This defines a homomorphism
$\Lambda
\to E$. We can extend this to a homomorphism $\widetilde \pi\colon \widetilde
\Lambda
\to E$ in the following way. Given $\alpha \in \widetilde \Lambda$, let $k = \deg
(L_\alpha|D)$. Define $\widetilde \pi (\alpha)$ to be the point corresponding to
the line bundle $(\varphi^*)^{-1}(L_\alpha)\otimes \scrO_E(-kp_0)$. Then
$\widetilde \pi$ is a homomorphism, and $\widetilde \pi(\kappa) = 0$. Hence
there is an induced homomorphism $\widetilde \Lambda/\Zee \cdot \kappa \to E$.
Since the intersection pairing identifies $(\widetilde \Lambda/\Zee \cdot
\kappa)^*$ with $\Lambda$, we see that this homomorphism is equivalent to an
element $\pi$ of $E\otimes_\Zee \Lambda$. Changing the choice of the marking
$\psi$ amounts to acting on $\Lambda$ via $W$, so that the period point is
invariantly defined in $(E\otimes_\Zee \Lambda)/W$.

Suppose that $\widetilde \xi$ is an $\widetilde {\mathbf{E}}_r$-bundle over $E$ and
that $\lambda$ is a
$\Cee^*$-bundle over
$E$. The inclusion of $\Cee^*$ into the center of $\widetilde {\mathbf{E}}_r$
allows us to define a new bundle $\widetilde \xi\otimes \lambda$. Let
$\det \widetilde \xi$ be the
$\Cee^*$-bundle induced by the surjection $\widetilde {\mathbf{E}}_r \to \Cee^*$.
Note that 
$$\det(\widetilde \xi\otimes \lambda) = \det \widetilde \xi \otimes
\lambda^{9-r}.$$

\begin{theorem} Given a triple $(X,D, \varphi)$, let $\varphi^*(\Xi|D)=\widetilde
\xi$. Then there is a natural choice of  a $\Cee^*$-bundle $\lambda$ for which
there is a
 canonical reduction of the structure group of the bundle 
$\widetilde \xi\otimes \lambda$ to an $\mathbf{E}_r$-bundle $\xi$ over
$E$. Moreover:
\begin{enumerate}
\item[\rm (i)] The bundle $\xi$ is regular, i.e. $h^0(E; \ad \xi) = r$.
\item[\rm (ii)]  The period point of the triple $(X,D, \varphi)$ in 
$(E\otimes_\Zee\Lambda)/W$ corresponds to the point defined by $\xi$, viewing 
$(E\otimes_\Zee
\Lambda)/W$ as the moduli space of semistable $\mathbf{E}_r$-bundles over $E$.
\end{enumerate}
\end{theorem}
\begin{proof}  It is easy to check via the exact sequence
$$\{1\} \to \mathbf{E}_r \to \widetilde {\mathbf{E}}_r \to \Cee^*\to \{1\}$$
that the obstruction to lifting $\widetilde \xi\otimes \lambda$
to a $\mathbf{E}_r$-bundle vanishes if and only if the bundle $\det(\widetilde
\xi\otimes
\lambda)$ is trivial. Moreover, using the above sequence
and the fact that there is a subgroup of $\Aut(\widetilde \xi\otimes \lambda)$
isomorphic to
$\Cee^*$ which surjects onto $\Aut(\det (\widetilde \xi\otimes \lambda))\cong
\Cee^*$   any lift of $\widetilde \xi\otimes
\lambda$ to
$\mathbf{E}_r$ is unique up to isomorphism. Also, it is easy to verify that $\det
(\widetilde \xi\otimes \lambda) = \det (\widetilde \xi_{\widetilde H}\otimes
\lambda)$, where $\widetilde \xi_{\widetilde H}$ is the corresponding
$\widetilde H$-bundle.

The exact sequence of tori
$$\{1\} \to H \to \widetilde H \to \Cee^* \to \{1\}$$
corresponds on the level of fundamental groups to 
$$0\to \Lambda \to \widetilde \Lambda \to \Zee\to 0,$$
where the surjection  $\widetilde \Lambda \to \Zee$ is given by $\lambda \mapsto
\langle \lambda, \kappa\rangle$. It follows from Lemma~\ref{taut} that $\det
\Xi_{\textrm{split}} = \det \Xi_{\widetilde H}$ corresponds to the line
bundle $(K_{\widetilde X})^{-1}$. Thus,
$\det \widetilde \xi$ corresponds to the line bundle $\scrO_E((9-r)p_0)$. It
follows that, if we set
$\lambda=\scrO_E(-p_0)$, then $\det
(\widetilde \xi\otimes \lambda)$ is trivial and thus $\widetilde \xi\otimes
\lambda$  lifts to a $\mathbf{E}_r$-bundle
$\xi$. The regularity of $\xi$ follows from the regularity of $\widetilde \xi$.

Finally, let us check that the period point of $(X,D, \varphi)$ corresponds to
the point it defines in the moduli space. Since this point is independent of the
S-equivalence class of $\xi$, we may as well work with
$\varphi^*(\Xi_{\textrm{split}}|D)\otimes \scrO_E(-p_0)$. Working instead with
the $\widetilde H$-bundle $\xi_{\widetilde H}$, note that a $\widetilde H$-bundle
over $E$ is the same as an element of $\Pic E\otimes \widetilde \Lambda$. The
bundle $\xi_{\widetilde H}$ is characterized by the property that, for
$\alpha\in \widetilde \Lambda \cong \widetilde\Lambda^*$, the $\Cee^*$-bundle
arising from the   homomorphism $\widetilde \Lambda \to \Zee$ is
$(\varphi^{-1})^*(L_\alpha|D)$. From this, it is clear that the period point
agrees with the moduli point.
\end{proof}

\section{Fundamental weights and the geometry of del Pezzo surfaces}

In this section, we assume that $4\leq r\leq 8$, so that the group $\mathbf{E}_r$
is simple. Also, we shall assume tacitly that $X=\widetilde X$, i.e. that we are
working with a smooth del Pezzo surface, unless otherwise specified.

 Fix the basis
$h, e_1,
\dots, e_r$ for
$\widetilde\Lambda$ and view elements of $\widetilde\Lambda$ as defining
homomorphisms $\Lambda \to
\Zee$ via the inner product (which we continue to assume is given by $h^2=1,
e_i^2=-1$). In this way, there is a homomorphism from $\widetilde \Lambda$ to the
set of weights on
$\Lambda$ and
$\lambda\in \widetilde \Lambda$ is the trivial weight if and only if $\lambda$ is
an integral multiple of $\kappa = 3h-\sum _ie_i$.

As before, we choose the simple coroots $\alpha_1 = e_1-e_2, \dots, \alpha _{r-1}
= e_{r-1}-e_r, \alpha _r = h-e_1-e_2-e_3$. By our conventions on the sign of 
intersection form, the root dual to a coroot $\alpha$ is confusingly identified
with $-\alpha \in \widetilde \Lambda$ via the intersection form. Let
$\{\varpi_i\}$ be the dual basis to
$\{\alpha_i\}$. For future reference, we give the list of all the positive coroots
and the highest root:

\begin{lemma} The positive roots are:
\begin{itemize} 
\item $e_i-e_j$ for $i<j$;
\item $h-e_i-e_j-e_k$;
\item $(r\geq 6)$ $2h- \sum _{i\in I}e_i$, where $\# I = 6$;
\item $(r= 8)$ $3h - 2e_i - \sum _{j\neq i}e_j$.
\end{itemize}

The highest root $\widetilde \alpha$ is given as follows: For $r=4,5$,
$\widetilde \alpha = h-e_{r-2} -e_{r-1} - e_r$, for $r=6,7$, $\widetilde \alpha
= 2h - \sum _{i=1}^6e_{r+1-i}$, for $r=8$,  $\widetilde \alpha
= 3h - \sum _{j=1}^7e_j - 2e_8 $.\qed
\end{lemma}

A routine calculation gives:

\begin{lemma} The dual basis element $\varpi_i$ is the image of the following
element $\widetilde\varpi_i$ of $\widetilde \Lambda$:
\begin{itemize}
\item $\widetilde\varpi_1 = h-e_1$, and the  linear system of all divisors on $X$
whose cohomology class is $h-e_1$ is a pencil of conics on
$X$;
\item $\widetilde\varpi_2= 2h-e_1-e_2$ and the  corresponding linear system is
defined by a smooth rational curve of degree $4$.
\item For $3\leq i\leq r-1$, $\widetilde\varpi_i = e_{i+1} + \cdots + e_r$
corresponds to
$r-i$ disjoint lines on $X$.
\item $\widetilde\varpi_r = h$ and the  corresponding linear system is defined by
a twisted cubic.
\qed
\end{itemize}
\end{lemma}

The above lemma shows that, for example, a lift of the fundamental weight
$\varpi_1$ may be identified with a particular pencil of conics on $X$. Of
course, $W$ acts on the set of all weights.  Via the marking, $W$ also acts on
$H^2(X; \Zee)$  and hence on the group of divisor classes of
$X$. In fact, $W$ acts simply transitively on the set of blowdowns $X\to \Pee^2$
together with a labeling of the exceptional curves. We next show that the Weyl
orbits of the corresponding fundamental weights may be identified with the
corresponding geometric objects on
$X$:

\begin{lemma} The Weyl group
$W$ operates transitively on any one of the following sets of objects: disjoint
sets of $k$ lines, $k\leq r-3$; pencils of conics; twisted cubics; linear systems
of smooth rational curves of degree four.
\end{lemma}
\begin{proof} Every set of $k\leq r-3$ disjoint lines  can be completed to a set
of $r$ disjoint lines. On the other hand, such a set defines a blowdown to
$\Pee^2$ and hence a diagonal basis $\gamma, \epsilon_1, \dots, \epsilon _r$ with
$\gamma^2=1, \epsilon _i^2=-1$ and $\kappa = 3\gamma -\sum _i\epsilon_i$. The
 Weyl group acts transitively on such bases, by Lemma~\ref{autom}, and hence on
the set of all sets of $k$ disjoint lines. A similar argument handles the case
of $r-1$ disjoint lines such that the blowdown is $\Pee^1\times \Pee^1$. 

Now suppose that $C$ is a twisted cubic on $X$, i.e. a smooth rational curve $C$
with $C^2=1$. The linear system $|C|$ defines a birational morphism to $\Pee^2$,
and we have seen that all such are conjugate under the Weyl group. If $C$ is a
conic, then it is easy to check that there is a blowdown to $\Pee^2$ such that
$C=\gamma-\epsilon_1$ for $\gamma$ the pullback of a hyperplane class and
$\epsilon_1$ an exceptional curve, and again all such are conjugate under $W$.
Finally, if $C$ is a smooth rational quartic, then $|C|$ defines a birational
morphism from $X$ to a quadric in $\Pee^3$, and since $X$ has no $-2$-curves, the
quadric is smooth. We conclude in this case by the last sentence of the preceding
paragraph.
\end{proof}

Each $\varpi_i$ is the highest weight of an irreducible representation of
$\mathbf{E}_r$, called a \textsl{fundamental representation}. We make some
comments on the fundamental representations.

\smallskip
\noindent \textbf{Minuscule representations.} For $r=4$, every
fundamental representation is minuscule. For $r=5$, the three minuscule
representations are the standard representation, with highest weight
$\varpi_1$, and the two half-spin representations,  with highest weights 
$\varpi_4$,
$\varpi_5$ respectively. For
$r=6$, the two minuscule
representations are  the two representations  with highest weights  $\varpi_1$,
$\varpi_5$. For $r=7$, the representation  with highest weight
$\varpi_6$ is the unique minuscule representation. For $r=8$, there are no
minuscule representations.

\smallskip
\noindent \textbf{The adjoint representation.} This is the unique
quasi-minuscule, non-minuscule representation. In all cases the highest weight is
the class of $\kappa - \widetilde \alpha$ (recall our conventions on the signs of
roots). For
$r=4$, the adjoint representation is not fundamental. For $r=5, 6, 7,8$, it is
respectively the representation with highest weight
$\varpi_2,\varpi_6, \varpi_1, \varpi_7$. 

\smallskip
\noindent \textbf{Dual representations.} Two irreducible representations are dual
if the negative of the lowest weight of one is equal to the highest weight of the
other. In particular, a representation is isomorphic to its dual if and only if
the lowest weight is the negative of the highest weight $\mu$, or equivalently in
terms of Weyl groups, if $-w_0(\mu) = \mu$, where 
$w_0\in W$ is the unique element sending positive roots to negative roots.
For
$\mu=\varpi_\alpha$, $\alpha$ a simple root, this is equivalent to:  $\alpha$ is
fixed by the outer automorphism
$-w_0$ of the root system. For $E_7$ and $E_8$, every outer automorphism is
trivial and this condition is automatic. In terms of the lattice $\widetilde
\Lambda$, suppose that we are given lifts of $\varpi_i$ and $\varpi_j$ to
elements $\widetilde\varpi_i$ and $\widetilde\varpi_j$ in $\widetilde \Lambda$.
Then the representations corresponding to
$\varpi_i$ and
$\varpi_j$ are dual if and only if the sum of
$\widetilde\varpi_i$ and a lift to $\widetilde \Lambda$ of the lowest weight
corresponding to
$\varpi_j$ add up to a multiple of $\kappa$. Clearly, this condition holds if and
only if there is some
$w\in W$ such that
$\widetilde\varpi_i + w(\widetilde\varpi_j) = n\kappa$. For example, in case
$r=6$, i.e. degree three, we can write the hyperplane section as the sum of a
line and a conic. For
$A_4$, the dualities amount to saying that we can write the hyperplane section
of a degree $5$ del Pezzo as the sum of a conic and a twisted cubic, or as a line
plus a smooth rational curve of degree four. For $D_5$, we can write the 
hyperplane section of a degree $4$ del Pezzo as the sum of a line and a twisted
cubic, as the sum of two conics, or as four lines, grouped in two pairs of
disjoint lines, and we can write twice the hyperplane section as a sum of two
smooth rational curves of degree four. Similar but more complicated results hold
for the smaller degrees.

\smallskip
\noindent \textbf{The cubic form.} Suppose that we are in the case $r=6$, i.e.
the case of cubic surfaces. The weights of the
$27$-dimensional representation correspond to lines on the cubic surface. As
was already known to Cartan \cite[pp.\ 272--273]{Cartan}, there is a symmetric
cubic form on the $27$-dimensional representation of $\mathbf{E}_6$, and in fact
$\mathbf{E}_6$ is the group of automorphisms of this form. In terms of the
geometry of the cubic surface, this form is determined by the following
condition: it is nontrivial on the tensor product of  three weight spaces    if
and only if the three corresponding lines sum up to a hyperplane section. It is an
exercise that the Weyl group acts transitively on the set of all sets of three
coplanar lines, so in fact up to a scalar multiple this is the only possible
invariant cubic form.

\smallskip

To give a partial answer to the question of Manin posed in the introduction, the
meaning of number of double sixes is that it is one-half the number of roots,
which we can see as follows. A double six consists of two disjoint
collections of lines $\{L_1, \dots, L_6\}$ and $\{L_1', \dots, L_6'\}$ such
that $L_i\cdot L_j = L_i'\cdot L_j'=0$ for all $i,j$, and $L_i\cdot
L_j'=\delta_{ij}$. Moreover, the six lines $\{L_1', \dots, L_6'\}$ are
determined by the set $\{L_1, \dots, L_6\}$. Instead of choosing a double six,
choose one of the two collections of disjoint lines that it determines---there are
$72$ of these. Each such determines a blowdown to
$\Pee^2$ and hence a root $2\gamma-\epsilon_1-\dots -\epsilon_6$, where the
$\epsilon_i$ are the exceptional curves of the blowdown and $\gamma$ is the
pullback of the hyperplane section. For example, for the standard blowdown we
get the root $2h-\sum _ie_i$. Since all roots are Weyl conjugate, we must get
them all this way. The remaining question, concerning the $45$ coplanar triples
of lines, is connected, perhaps less explicitly, to the cubic form on the standard
representation of
$\mathbf{E}_6$.
 
We now consider the vector bundle induced from $\Xi$ under a representation of
$\widetilde {\mathbf{E}}_r$. Because we assume that $X$ is smooth, $\Xi =
\Xi_{\textrm{split}}$. Let
$\rho\colon \mathbf{E}_r\to GL(N)$ be an irreducible representation, and suppose
that we have extended $\rho$ to a representation $\widetilde \rho\colon
\widetilde {\mathbf{E}}_r\to GL(N)$. Such an extension is equivalent to choosing an
integer
$b$ such that
$$\rho(c) = \exp(2\pi \sqrt{-1}b/d)\cdot \Id,$$
where $d=9-r$ is the order of $c$ (cf.\ (i) of Lemma~\ref{duals}).
For example, there is a unique choice of $b$ with $0\leq b\leq d-1$. The map 
$$(g,t) \in \mathbf{E}_r\times \Cee^* \mapsto \rho(g)\cdot t^b\Id$$
then gives an extension $\widetilde \rho$. It is clear from the example at the
end of Section 1 that, if $\rho$ is the irreducible with highest weight
$\varpi_{r-1}$ and $c$ is the corresponding central element, then $\rho(c) =
\exp(2\pi \sqrt{-1}/d)\cdot \Id$ and in fact $\widetilde {\mathbf{E}}_r$ is just
$\rho(\mathbf{E}_r)\cdot \Cee^* \subseteq GL(N)$. In this case, we can always
choose the extension $\widetilde \rho$ to be the inclusion. More generally, using
(ii) of Lemma~\ref{duals}, an extension of $\rho$ to $\widetilde \rho$
corresponds to choosing a
$\lambda \in
\widetilde \Lambda$ such that the highest weight of $\rho$ is the image of
$\lambda$ under the map from $\widetilde \Lambda$ to $\Lambda^*$ given by the
inner product, and any two such differ by a multiple of $\kappa$.

For a
weight
$\lambda$ of $\widetilde\rho$, which we view via the intersection form as an
element of $\widetilde\Lambda$, the corresponding line bundle is the
line bundle we have denoted by $L_\lambda$. We have:

\begin{proposition} Suppose that $X$ is smooth. Let $\widetilde \rho\colon
\widetilde {\mathbf{E}}_r
\to GL(N)$ be the representation  defined by the lift $e_r$ of $\varpi_{r-1}$.
Then, for
$r\neq 8$, the induced vector bundle $\Xi\times _{\widetilde {\mathbf{E}}_r}\Cee^N$
is
$\bigoplus _i\scrO_X(L_i)$, where the $L_i$ are the distinct lines on $X$. For
$r=8$, $\Xi\times _{\widetilde {\mathbf{E}}_8}\Cee^{248} \cong \bigoplus
_i\scrO_X(L_i)\oplus (K_X^{-1})^{\oplus 8}$.
\end{proposition}
\begin{proof} For $r\neq 8$, this is clear by the remarks before the statement.
For $r=8$, the adjoint bundle $\ad \Xi \cong \bigoplus_{\alpha\in R}L_\alpha
\oplus \scrO_X^8$, where $R\subseteq H^2(X;\Zee)$ is the set of roots. The result
then follows since
$\Xi\times _{\widetilde {\mathbf{E}}_8}\Cee^{248} \cong \ad \Xi \otimes K_X^{-1}$.
\end{proof}

For $r=6$, there is also the natural cubic form on the rank $27$ vector bundle,
with values in the line bundle $\scrO_X(D) =K_X^{-1}$, where $D$ is a hyperplane
section of $X$, defined by the obvious map
$$\scrO_X(L_i) \otimes \scrO_X(L_j) \otimes \scrO_X(L_k) \to \scrO_X(L_i+
L_j+L_k).$$ 

Finally, we discuss the structure of the bundle in case $X$ has rational double
points. For simplicity, let us just consider the case where $r=6$ and there is
exactly one double point whose preimage on $\widetilde X$ is $C$. For example, we
could take
$[C]=e_1-e_2$. The $27$ classes $e_i, h-e_i-e_j, 2h-\sum _{i\neq k}e_i$ which
define lines on a smooth cubic---call these \textsl{numerical lines}---specialize
to $6$ pairs of the form $\{\lambda, \lambda +C\}$ corresponding to $6$ lines on
the singular cubic $X$ which pass through the double point, together with $15$
lines which do not meet the double point. Let $L_i'$, $i=1, \dots, 6$ be the line
bundles corresponding to numerical lines $\lambda$ such that $\lambda +C$ is also
a numerical line. This is equivalent to: $L_i'\cdot C = 1$. Let $L_j''$, $j=1,
\dots, 15$ be the remaining classes, so that $L_j''\cdot C =0$. Finally, let $V_2$
be the unique rank two vector bundle on $\widetilde X$ which is a nontrivial
extension
$$0 \to \scrO_{\widetilde X}(C) \to V_2 \to \scrO_{\widetilde X} \to 0.$$
It is easy to check that the map $H^1(\widetilde X;\scrO_{\widetilde X}(C)) \to
H^1(C; \scrO_{\widetilde X}(C)|C) = H^1(C; \scrO_C(-2))$ is an isomorphism. Thus,
the extension $V_2$ restricts on $C$ to the nontrivial extension
$$0 \to \scrO_C(-2) \to V_2|C \to \scrO_C\to 0.$$
It follows that $V_2|C \cong \scrO_C(-1) \oplus \scrO_C(-1)$. Thus, if we consider
the rank $27$ bundle
$$\mathcal{V}=\bigoplus _{i=1}^6(V_2\otimes L_i') \oplus \bigoplus
_{j=1}^{15}L_j'',$$
then $\mathcal{V}$ restricts to the trivial rank $27$ bundle over $C$ and thus is
pulled back from $X$. It is easy to see that $\mathcal{V}$ is the bundle induced
from the standard representation of $\widetilde {\mathbf{E}}_6$. Similar but
considerably more elaborate descriptions exist for the case of arbitrary rational
double points.

\bigskip
\noindent
Department of Mathematics \\
Columbia University \\
New York, NY 10027 \\
USA

\bigskip
\noindent
{\tt rf@math.columbia.edu, jm@math.columbia.edu}

\end{document}